\documentclass[12pt]{article}
\usepackage{amsmath}
\usepackage{amssymb}
\usepackage{amscd}
\usepackage{amsthm}

\theoremstyle{plain}
\newtheorem{Thm}{Theorem}[section]

\newtheorem{Cor}[Thm]{Corollary}
\newtheorem{Lem}[Thm]{Lemma}

\theoremstyle{definition}

\newtheorem{Expl}[Thm]{Example}

\numberwithin{equation}{section}

\title{Derived Categories and Birational Geometry}
\author{Yujiro Kawamata}

\begin{document}

\maketitle

\section{Introduction}
This paper is concerned with a surprizing parallelism 
between minimal model program and 
semi-orthogonal decompositions of derived categories found by 
Bondal and Orlov (\cite{BO}).

The bounded derived category of coherent sheaves on a smooth projective 
variety has a finiteness property called saturatedness (\S 2).
A map between smooth projective varieties 
such as a projective space bundle, a blowing up along a smooth center or 
a standard smooth flip 
induces a semi-orthogonal decomposition of the derived category (\S 3).
We expect that one can attach a saturated triangulated category 
to any singular variety which appears in the minimal model program,
and that a map between these varieties such as a Mori fiber space, 
a divisorial 
contraction or a flip induces a semi-orthogonal decomposition
of the category.

We review some basic results on categories in \S 2 following \cite{BK}.
We treat some cases which motivate our expectation in \S 3.
The cases include smooth varieties due to \cite{Orlov} and \cite{BO}, 
and toric varieties \cite{toric}.
We provide some properties of the desired category, and 
give a possible definition of the correct category for certain three 
dimensional cases by calculating the behavior of the category under the  
divisorial contractions of smooth threefolds in \S 4.
We conclude this paper with a short section of open questions.

The author would like to thank the referee for a careful reading and useful 
suggestions.


\section{Saturated categories and semi-orthogonal decompositions}

Let $X$ be a projective variety over a field $k$.
The largest derived category for $X$ is the unbounded derived 
category of 
quasi-coherent sheaves denoted by $D(\text{QCoh}(X))$ 
(we refer to \cite{GM} for basic definitions).
There are functors such as $f_*,f^*,f^!$ for morphisms $f: X \to Y$, 
and the Grothendieck duality theorem holds (\cite{Neeman}).
This category is {\em $k$-linear} in the sense that the 
set of morphisms $\text{Hom}(a,b)$ for $a,b \in D(\text{QCoh}(X))$ 
has the structure of a $k$-vector space.
It is, however, infinite dimensional and we might hope to work with 
categories for which the Hom space has finite dimension.

A candidate is the bounded derived category of coherent sheaves denoted 
by $D^b(\text{Coh}(X))$.
If $X$ is smooth, then $D^b(\text{Coh}(X))$ is of {\em finite type} 
in the sense that 
$\sum_{p \in \mathbb{Z}} \dim \text{Hom}(a,b[p])$ is finite.
But if $X$ is singular, then the homological dimension is infinite. 
Indeed, 
if $x \in X$ is a singular point, then there are infinitely many $p$ 
such 
that $\text{Hom}(\mathcal{O}_x, \mathcal{O}_x[p]) \ne 0$, where 
$\mathcal{O}_x$ denotes a skyscraper sheaf of length $1$ at $x$.

An object $c$ in a triangulated category $D$ having arbitrary coproducts 
is said to be {\em compact} if 
$\text{Hom}(c, \coprod_{\lambda} a_{\lambda})
\cong \coprod_{\lambda}\text{Hom}(c, a_{\lambda})$.
An object $c \in D(\text{QCoh}(X))$ is compact if and only if 
it is {\em perfect} in the sense that it is locally isomorphic to 
a bounded complex of locally free coherent sheaves.
Let $\text{Perf}(X) = D(\text{QCoh}(X))^c$ be the triangulated 
category of perfect complexes.
It is a full subcategory of $D^b(\text{Coh}(X))$, and 
they coincide if $X$ is smooth.

If $a \in \text{Perf}(X)$ and $b \in D^b(\text{Coh}(X))$, then
$\sum_{p \in \mathbb{Z}} \dim \text{Hom}(a,b[p]) < \infty$.
Serre duality holds
\[
\text{Hom}(a,b) \cong \text{Hom}(b,a \otimes \omega^{\bullet}_X)^*
\]
for $a \in \text{Perf}(X)$ and $b \in D^b(\text{Coh}(X))$,
where $\omega^{\bullet}_X$ is a dualizing complex.
For example, if $X$ is smooth, then $\omega^{\bullet}_X = 
\omega_X[\dim X]$.
Therefore, $\text{Perf}(X)$ is similar to homology while 
$D^b(\text{Coh}(X))$
cohomology, because Serre duality is similar to Poincar\'e 
duality.
In particular, if $X$ is smooth, then the functor
$S_X: D^b(\text{Coh}(X)) \to D^b(\text{Coh}(X))$ defined by 
$S_X(a) = a \otimes \omega_X[\dim X]$ is a {\em Serre functor}
in the sense that 
\begin{equation}\label{Serre}
\text{Hom}(a,b) \cong \text{Hom}(b,S_X(a))^*
\end{equation}
for $a,b \in D^b(\text{Coh}(X))$.

It is important to consider singular varieties in minimal model
theory.
Since the sets $\text{Hom}(\mathcal{O}_x, \mathcal{O}_x[p])$ 
do not vanish for 
$p \ge 0$ but vanish for $p < 0$ for a singular point $x \in X$, 
there cannot be a Serre functor in the category $D^b(\text{Coh}(X))$
for a singular variety $X$.
So we would like to have something like intersection homology 
for a singular variety $X$ that lies between
$\text{Perf}(X)$ and $D^b(\text{Coh}(X))$ and has a Serre functor.

A $k$-linear triangulated category $D$ of finite type is said to be 
{\em saturated} 
if any cohomological functor $F: D^{\text{op}} \to (k\text{-vect})$ 
of finite type is representable.
That is, if $\sum_{p \in \mathbb{Z}} \dim F(a[p]) < \infty$ for any
$a \in D$, then there exists $b \in D$ such that 
$F(a) \cong \text{Hom}(a,b)$.
A saturated category has always a Serre functor, because the dual of 
the left hand side of (\ref{Serre}) defines a cohomological functor.
The category $D^b(\text{Coh}(X))$ for smooth $X$ is saturated 
(\cite{BvdB}).
For example, 
an object $a \in D$ in a $k$-linear triangulated category is said to be 
{\em exceptional}
if $\text{Hom}(a,a[p]) = 0$ for $p \ne 0$ while 
$\text{Hom}(a,a) \cong k$.
In this case, the subcategory $\langle a \rangle$ generated by $a$ is 
saturated because it is equivalent to $D^b(\text{Coh}(\text{Spec }k))$.

A triangulated full subcategory $B$ of a triangulated category $A$ 
is said to be {\em right admissible} if the embedding $j_*: B \to A$  
has a right adjoint functor $j^!: A \to B$.
Let $C = B^{\perp}$ be the {\em right orthogonal} defined by
\[
C = \{c \in A \,\, \vert \,\, \text{Hom}(b,c) = 0 \,\, 
\forall b \in B\}.
\]
Then an arbitrary $a \in A$ has a unique presentation by a 
distinguished triangle
\[
b \to a \to c \to b[1]
\]
for $b = j_*j^!a \in B$ and $c \in C$.
We write this as
\[
A = \langle C, B \rangle
\]
and call it a {\em semi-orthogonal decomposition} of $A$.
We note that the form $\text{Hom}(a,b)$ is similar to a bilinear 
form on $A$, but it is not symmetric.
Thus the orthogonality is only one-sided.
Indeed, $D^b(\text{Coh}(X))$ is indecomposable with respect to 
$\text{Hom}$ as long as $X$ is 
irreducible (\cite{Bridgeland}).
We denote by $\langle D, C, B \rangle$ for 
$\langle D,\langle C, B\rangle \rangle$, etc.
For example, if $A$ is of finite type and $B$ is saturated, 
then $B$ is always right admissible.
Conversely, if $A$ is saturated and $B$ is right admissible, 
then $B$ is saturated.

If $A$ has a Serre functor $S_A$, then 
$B$ has also a Serre functor given by $S_B = j^!S_Aj_*$,
because
\[
\begin{split}
&\text{Hom}_B(a,b) \cong \text{Hom}_A(j_*a,j_*b) \\
&\cong \text{Hom}_A(j_*b,S_Aj_*a)^*
\cong \text{Hom}_B(b,j^!S_Aj_*a)^*
\end{split}
\]
for $a,b \in B$.
$j_*$ has also a left adjoint $j^*$ defined by 
$j^* = S_B^{-1}j^!S_A$, because
\[
\begin{split}
&\text{Hom}_A(a,j_*b) \cong \text{Hom}_A(j_*b,S_Aa)^* \\
&\cong \text{Hom}_B(b,j^!S_Aa)^*
\cong \text{Hom}_B(j^!S_Aa,S_Bb)
\end{split}
\]
for $a \in A$ and $b \in B$.

\begin{Lem}\label{saturated}
Assume that $X$ is singular.
Then $\text{Perf}(X)$ is not saturated.
\end{Lem}

\begin{proof}
Assume the contrary and let $x \in X$ be a singular point. 
Assume that the contravariant functor 
$\text{Hom}(\bullet, \mathcal{O}_x)$ 
on $\text{Perf}(X)$ is represented by an object $c \in \text{Perf}(X)$:
\[
\text{Hom}(\bullet, \mathcal{O}_x) \cong \text{Hom}(\bullet,c).
\]
Let $x' \ne x$ be another point and $b$ a coherent sheaf 
supported at $x'$ such that $b \in \text{Perf}(X)$ as an object in 
$D^b(\text{Coh}(X))$. 
Then we have 
\[
\text{Hom}(b, c) \cong \text{Hom}(b, \mathcal{O}_x) = 0.
\]
Hence $\text{Supp}(c) = \{x\}$.
Since 
\[
\text{Hom}(\mathcal{O}_X, c) 
\cong \text{Hom}(\mathcal{O}_X, \mathcal{O}_x) 
\cong k
\]
we conclude that 
$c$ is a sheaf of length $1$ supported at $x$, 
but the latter is not contained in $\text{Perf}(X)$.
\end{proof}


\section{Minimal model program and decompositions}

The minimal model program consists of three kinds of basic operations, 
namely, Mori fiber spaces, divisorial contractions and flips
(see for example \cite{KMM}).
The latter two are birational maps which decrease 
the canonical divisor $K$.
The following are the simplest examples of these operations 
for smooth projective varieties.
The point is that semi-orthogonal decompositions of 
the derived categories
are parallel to the decompositions of canonical divisors.

\begin{Expl}
In this example, we consider only smooth projective 
varieties $X$,$Y$, etc.
For simplicity, we write $D(X)$ for the bounded derived category of 
coherent sheaves $D^b(\text{Coh}(X))$, etc.

(1) (\cite{Orlov}) Let $f: X \to Y$ be a projective space bundle 
associated to a vector bundle of rank $r$.
This is a Mori fiber space.
The functor $f^*: D(Y) \to D(X)$ is fully faithful, 
and there is a semi-orthogonal decomposition
\[
D(X) = \langle D(Y)_{-r+1}, \dots, D(Y)_{-1}, D(Y)_0 \rangle
\]
where the $D(Y)_i$ denote the subcategories 
$f^*D(Y) \otimes \mathcal{O}_X(i)$ 
of $D(X)$ for the tautological line bundle $\mathcal{O}_X(1)$.

(2) (\cite{Orlov}) Let $C$ be a smooth subvariety of $Y$ of codimension 
$r \ge 2$
and $f: X \to Y$ the blowing up along $C$.
This is a divisorial contraction.
The exceptional divisor $E$ is a projective space bundle 
over $C$ as in (1).
Let $f_E: E \to C$ be the induced morphism and $j: E \to X$ 
the embedding.
Then the functors $f^*: D(Y) \to D(X)$ and $j_*f_E^*: D(C) \to D(X)$ 
are 
fully faithful, and there is a semi-orthogonal decomposition
\[
D(X) = \langle j_*f_E^*D(C)_{-r+1}, \dots, j_*f_E^*D(C)_{-1}, 
f^*D(Y) \rangle
\]
where $j_*f_E^*D(C)_i = j_*f_E^*D(C) \otimes \mathcal{O}_X(-iE)$.
We have a corresponding equality of canonical divisors 
\[
K_X = f^*K_Y + (r-1)E.
\]

(3) (\cite{BO}) Let $E$ be a subvariety of $X$ which is isomorphic to 
a projective space $\mathbb{P}^{r-1}$ 
and such that the normal bundle $N_{E/X}$ is 
isomorphic to $\mathcal{O}_{\mathbb{P}^{r-1}}(-1)^s$.
Let $f: Y \to X$ be the blowing up along $E$.
Then there is a blowing down $f^+:Y \to X^+$ of the exceptional divisor 
$F = f^{-1}(E)$ to another direction such that $E^+ = f^+(F)$ is 
isomorphic to $\mathbb{P}^{s-1}$ and such that the normal bundle 
$N_{E^+/X^+}$ 
is isomorphic to $\mathcal{O}_{\mathbb{P}^{s-1}}(-1)^r$.
If $r > s$, then this is a flip, while if $r=s$, then it is a flop.
If $r \ge s$, then the functor $f_*f^{+*}: D(X^+) \to D(X)$ is 
fully faithful, and there is a semi-orthogonal decomposition
\[
D(X) = \langle \mathcal{O}_E(s-r), \dots, \mathcal{O}_E(-1), 
f_*f^{+*}D(X^+) \rangle
\]
where the subcategories 
$\langle \mathcal{O}_E(i) \rangle$ generated by the sheaves 
$\mathcal{O}_E(i)$ are denoted by $\mathcal{O}_E(i)$ for simplicity.  
In particular, we have an equivalence of triangulated categories
$D(X) \cong D(X^+)$ for the flop case.
We have a corresponding equality of canonical divisors 
\[
f^*K_X = f^{+*}K_{X^+} + (r-s)F.
\]
\end{Expl}

It is important to deal with singular varieties 
in the minimal model program.
Therefore, we have to define good derived categories for such
varieties.
The simplest case is the one with quotient singularities.
For a variety $X$ with only quotient singularities, we can naturally associate
a smooth Deligne-Mumford stack $\mathcal{X}$.
The set of points of the stack $\mathcal{X}$ is the same as that of the 
variety $X$, but
the points on $\mathcal{X}$ have automorphism groups 
corresponding to
the stabilizer groups of the points on $X$.
The sheaves on $\mathcal{X}$ have actions by these groups.
Let $D^b(\text{Coh}(\mathcal{X}))$ be the bounded derived category of 
coherent sheaves on $\mathcal{X}$.
For example, if $X = M/G$ is the quotient of a smooth variety $M$ by a 
finite group $G$, 
then $D^b(\text{Coh}(\mathcal{X})) = D^b(\text{Coh}^G(M))$ is 
the bounded derived category of $G$-equivariant coherent sheaves on $M$.
The following example suggests that the above 
$D^b(\text{Coh}(\mathcal{X}))$ is the
correct answer to our problem for varieties 
with only quotient singularities (cf. \cite{log-crep} and \cite{toric} 
for more justifications).

\begin{Expl} (\cite{Francia})
Let $X$ be a $4$ dimensional smooth projective variety, 
$E$ a subvariety 
which is isomorphic to a projective plane 
$\mathbb{P}^2$ and such that the normal bundle 
$N_{E/X}$ is isomorphic to $\mathcal{O}_{\mathbb{P}^2}(-1) \oplus
\mathcal{O}_{\mathbb{P}^2}(-2)$.
Let $f_1: X_1 \to X$ be the blowing up along $E$.
Then the exceptional divisor $F = f_1^{-1}(E)$ contains a subvariety 
$E_1$ which is isomorphic to $\mathbb{P}^2$ and the normal bundle 
$N_{E_1/X_1}$ is isomorphic to $\mathcal{O}_{\mathbb{P}^2}(-1)^2$.
Let $f_2:Y \to X_1$ be the further blowing up along $E_1$.
Then there is a blowing down $f_1^+:Y \to X^+_1$ 
of the exceptional divisor 
$F_1 = f_2^{-1}(E_1)$ to another direction such that 
$E_1^+ = f_1^+(F_1)$ is 
isomorphic to $\mathbb{P}^1$.
The strict transform $F'$ of $F$ on $X^+_1$ is isomorphic to 
$\mathbb{P}^2$
and the normal bundle $N_{F'/X^+_1}$ is isomorphic to 
$\mathcal{O}_{\mathbb{P}^2}(-2)$.
Let $f^+_2: X^+_1 \to X^+$ be the blowing down of $F'$.
The image $Q = f^+_2(F')$ is an isolated quotient singularity 
with stabilizer group $\mathbb{Z}/2$.

The composition $f_2^+f_1^+f_2^{-1}f_1^{-1}$ is a flop.
Namely, we have an equality $f_2^*f_1^*K_X = f_1^{+*}f_2^{+*}K_{X^+}$.
Correspondingly, we have an equivalence of derived categories
\[
f_{2*}^+f_{1*}^+f_2^*f_1^*: 
D^b(\text{Coh}(X)) \cong D^b(\text{Coh}(\mathcal{X}^+))
\]
where we consider the associated stack $\mathcal{X}^+$ 
instead of the underlying variety with a quotient singularity $X^+$.
But the functor $\pi_*: D^b(\text{Coh}(\mathcal{X}^+)) 
\to D^b(\text{Coh}(X^+))$
induced by the projection $\pi: \mathcal{X}^+ \to X^+$ 
is not an equivalence.
For example, 
\[
f_{2*}^+f_{1*}^+f_2^*f_1^*(\Omega^1_E(-1)) = \mathcal{O}_Q(1)
\]
and $\pi_*\mathcal{O}_Q(1) \cong 0$, where $\mathcal{O}_Q(1)$ 
is a skyscraper
sheaf of length $1$ supported at $Q$ on which the stabilizer 
group acts non-trivially.
Thus an ordinary sheaf $\Omega^1_E(-1)$ on $X$ corresponds to an equivariant
sheaf $\mathcal{O}_Q(1)$ on the imaginary stack $\mathcal{X}^+$ 
which disappears on the real variety $X^+$.
\end{Expl}

The following example shows that the semi-orthogonal decompositions of 
derived categories are governed by the inequalities of canonical 
divisors and not by the directions of morphisms.
This fact suggests a distinguished status of the canonical divisors 
in the theory of derived categories.
We note that the derived categories contain almost all information on 
varieties in a similar way like the motives (cf. \cite{Orlov2}).

\begin{Expl} (\cite{toric})
Let $X$ be a smooth projective variety of dimension $n$ 
which contains a divisor $E$ being isomorphic to a projective space 
$\mathbb{P}^{n-1}$ and such that the normal bundle is isomorphic to 
$\mathcal{O}_{\mathbb{P}^{n-1}}(- k)$ for an integer $k > 0$.
Let $f: X \to Y$ be the blowing down of $E$.
Then $Y$ has an isolated quotient singularity $Q$ 
whose stabilizer group is 
isomorphic to $\mathbb{Z}/k$.
Let $\mathcal{Y}$ be the associated smooth Deligne-Mumford stack,
and let $\mathcal{Z} = X \times_Y \mathcal{Y}$ be the fiber product with
projections $\pi: \mathcal{Z} \to X$ and 
$\tilde f: \mathcal{Z} \to \mathcal{Y}$.
We have an equality
\[
K_X = f^*K_Y + \frac {n-k}kE.
\]
If $n > k$, then we have $K_X > f^*K_Y$.
Correspondingly, the functor $\pi_*\tilde f^*: 
D^b(\text{Coh}(\mathcal{Y}))
\to D^b(\text{Coh}(X))$ is fully faithful, 
and there is a semi-orthogonal decomposition
\[
D^b(\text{Coh}(X)) = \langle \mathcal{O}_E(-n+k), 
\dots, \mathcal{O}_E(-1), 
\pi_*\tilde f^*D^b(\text{Coh}(\mathcal{Y})) \rangle.
\]
On the other hand, if $n < k$, then we have $K_X < f^*K_Y$.
Correspondingly, the functor $\tilde f_*\pi^*: D^b(\text{Coh}(X))
\to D^b(\text{Coh}(\mathcal{Y}))$ is fully faithful, 
and there is a semi-orthogonal decomposition
\[
D^b(\text{Coh}(\mathcal{Y})) = \langle \mathcal{O}_Q(-n), \dots, 
\mathcal{O}_Q(-k+1),\tilde f_*\pi^*D^b(\text{Coh}(X)) \rangle
\]
where $\mathcal{O}_Q(i)$ denotes a skyscraper sheaf of length $1$ at $Q$ 
with a suitable action by the stabilizer group.
In particular, if $n=k$, then there is an equivalence
$\pi_*\tilde f^*: D^b(\text{Coh}(\mathcal{Y}))
\cong D^b(\text{Coh}(X))$.

We note that the functors are given by the pull-backs 
and the push-downs as in the case of flips.
Indeed, divisorial contractions and flips are very similar operations 
from the view point of the minimal model program.
\end{Expl}

The above picture extends for $\mathbb{Q}$-factorial 
toric varieties (\cite{toric}):

\begin{Thm}
Let $f: X -\to Y$ be a toric divisorial contraction or flip between 
$\mathbb{Q}$-factorial projective toric varieties, and
let $\mathcal{X}$ and $\mathcal{Y}$ be their associated 
smooth Deligne-Mumford stacks.
Then there is a semi-orthogonal decomposition
\[
D^b(\text{Coh}(\mathcal{X})) \cong 
\langle C, D^b(\text{Coh}(\mathcal{Y})) \rangle
\]
which is described in detail in terms of toric fans for $X$ and $Y$.
\end{Thm}

We note that $\mathbb{Q}$-factorial toric varieties have only 
quotient singularities.
We have similar results for toric Mori fiber spaces.
The theorem has a log version in the case where 
the coefficients of the boundary are of the form $1 - \frac 1m$ 
for some positive integers $m$.
In particular, toric flops induce derived equivalences.
We refer to \cite{toric} for details.
As a corollary, we obtain the McKay correspondence for abelian 
quotient singularities:

\begin{Cor}
Let $X$ be a projective variety with only quotient singularities 
whose stabilizer groups are abelian groups 
whose orders are prime to the characteristic of the base field, 
and let $f: Y \to X$ be a projective crepant resolution, i.e., 
$Y$ is smooth, $f$ is projective and birational, and $K_Y = f^*K_X$.
Then there is an equivalence of triangulated categories
\[
D^b(\text{Coh}(Y)) \cong D^b(\text{Coh}(\mathcal{X})).
\]
\end{Cor}

\begin{proof}
We may replace $X$ by its local model and assume that $X$ is toric by 
\cite{Bridgeland}.
Then there exists a toric crepant $\mathbb{Q}$-factorial terminalization 
$f': Y' \to X$.
Let $H$ be an $f$-ample divisor on $Y$ and $H'$ its strict transform on $Y'$.
Since $Y'$ is toric, $H'$ is linearly equivalent to a toric divisor 
which is denoted by $H'$ again.
We proceed by MMP with respect to $(Y', \epsilon H')$ over $X$ 
for a small positive number $\epsilon$.
Since the pair is toric, the process is a toric MMP.
After finitely many steps, we reach a log minimal model
that is isomorphic to $Y$.
Therefore, $f: Y \to X$ is also toric.
The MMP over $X$ starting from $(Y, B)$ for a 
suitable toric boundary $B$ ends at $X$.
Therefore, $Y$ and the smooth Deligne-Mumford stack over $X$ are 
derived equivalent by the theorem (see also \cite{log-crep}).
\end{proof}

There is a warning.
The derived category may have semi-orthogonal decompositions 
beyond the minimal model program.
For example, the derived category of a projective space 
has a complete semi-orthogonal decomposition 
to exceptional objects by \cite{Beilinson}.
This fact is extended to an arbitrary $\mathbb{Q}$-factorial projective 
toric variety (\cite{toric}).
The derived categories of some Fano manifolds have interesting 
semi-orthogonal decompositions which reflect the geometry of these
manifolds (\cite{Kuznetsov}).
Even minimal varieties such as Enriques surfaces have derived
semi-orthogonal decompositions.
Therefore, MMP is only a preparation as in the case of the classification 
of surfaces.
After that, deeper decompositions may be possible like in the case of 
decompositions of motives.


\section{Divisorial contractions of smooth $3$-folds}

We would like to define a correct category for an arbitrary variety 
which 
appears in the minimal model program, or in the Mori category.
Thus let $X$ be a projective variety with only terminal singularities, and 
$f: Y \to X$ a resolution of singularities.
As a working hypothesis, 
we look for a minimal saturated subcategory $D = D(X)$ of $D^b(\text{Coh}(Y))$ 
which contains $f^*\text{Perf}(X)$:
\[
f^*\text{Perf}(X) \subset D(X) \subset D^b(\text{Coh}(Y))
\]
if it exists and unique up to equivalence.
If $X$ is smooth, then we have $D(X) = D^b(\text{Coh}(X))$.

As a first step, we consider a divisorial contraction $f: Y \to X$ 
of a smooth $3$ dimensional variety to a singular variety.
The morphism $f$ is an isomorphism outside a prime divisor $E$ on $Y$, 
and $P=f(E)$ is the singular point of $X$.
There are three cases, where two of them has already answers described in 
\S 3.
The above working hypothesis seems to produce the same categories 
in these cases.

\subsection{Case 1}

$E$ is isomorphic to a smooth quadric surface 
$\mathbb{P}^1 \times \mathbb{P}^1$ with normal bundle 
$\mathcal{O}_E(-1,-1)$, where the pair of rulings on $E$ 
are numerically equivalent to each other on $X$.

Let 
\[
C = \langle \mathcal{O}_E(-1,-1), \mathcal{O}_E(0,-1) \rangle 
\subset D^b(\text{Coh}(Y))
\]
be the subcategory generated by a sequence of exceptional objects,
and let $D = {}^{\perp}C$.
Since $C$ is equivalent to $D^b(\text{Coh}(\mathbb{P}^1))$,
it is admissible, and $D$ is saturated.
Since $f_*c = 0$ for any $c \in C$, 
we have $f^*p \in D$ for $p \in \text{Perf}(X)$, thus
$f^*\text{Perf}(X) \subset D$.

If $k = \mathbb{C}$, then there are two small resolutions 
$g_i: Y_i \to X$ ($i=1,2$) in the analytic category, 
with analytic divisorial contractions $f_i: Y \to Y_i$ 
such that $f_i(E)$ is isomorphic to $\mathbb{P}^1$.
We can check that $D = f_1^*D^b(\text{Coh}(Y_1))$, and the latter should be
the correct category since $Y_1$ is smooth.
Since $g_1$ is crepant, it follows that $S_D(d) \cong d[3]$ if 
$d \in D$ and $f_*d=0$.
We can also prove this fact for general $k$, 
because these objects are concentrated on the divisor $E$ 
and the global structure of $X$ is irrelevant.

\begin{Lem}
$D$ is minimal in the sense that 
$D$ has no semi-orthogonal decomposition relative to $X$, i.e., a 
semi-orthogonal decomposition such that 
one of the factors contains $f^*\text{Perf}(X)$.
\end{Lem}

\begin{proof}
Suppose there is still a semi-orthgonal decomposition
$D = \langle C',D' \rangle$ such that $f^*\text{Perf}(X) \subset D'$.
Since $C' \subset (f^*\text{Perf}(X))^{\perp}$, we have $f_*c = 0$ for 
$c \in C'$.
Then we have for $c \in C'$ and $d \in D'$
\[
\text{Hom}(c,d) \cong \text{Hom}(d,c[3])^* \cong 0.
\]
Hence $D$ is decomposable; $D = C' \oplus D'$.

Assuming that $k = \mathbb{C}$, 
let $a_y = f_1^*\mathcal{O}_y \in D$ for a point 
$y \in Y_1$.
Since $\text{Hom}(a_y,a_y) \cong k$, we have either $a_y \in C'$ or 
$a_y \in D'$.
Since $f_*a_y \ne 0$, we conclude that $a_y \in D'$.
But let $c' \in C'$ be a non-zero object, and write $c' = f_1^*c$.
If we take a point $y$ in the support of $c$, then we have 
$\text{Hom}(c,\mathcal{O}_y[p]) \ne 0$ for some $p$,
hence $\text{Hom}(c',a_y[p]) \ne 0$.
But this is a contradiction.
For general $k$, we can still define $a_y$ 
because it is supported on $E$,
and the above argument works.
\end{proof}

If we take 
\[
C_1 = \langle \mathcal{O}_E(-1,-1), \mathcal{O}_E(-1,0) \rangle 
\subset D^b(\text{Coh}(Y))
\]
instead of $C$ and define $D_1 = {}^{\perp}C_1$, then 
$D$ and $D_1$ are equivalent, because $Y_1$ and $Y_2$ are related by a 
standard flop.

There is another candidate for an admissible subcategory 
which is more symmetric than $C$ or $C_1$.
Let 
\[
\tilde C = \langle \mathcal{O}_E(-1,0), \mathcal{O}_E(0,-1) \rangle 
\subset D^b(\text{Coh}(Y))
\]
and let $\tilde D = {}^{\perp}\tilde C$ be the right orthogonal.
We have again $f^*\text{Perf}(X) \subset \tilde D$.
The generators of $\tilde C$ are exceptional objects,
\[
\text{Hom}(\mathcal{O}_E(-1,0), \mathcal{O}_E(0,-1)[2]) 
\cong \text{Hom}(\mathcal{O}_E(0,-1), \mathcal{O}_E(-1,0)[2]) 
\cong k
\]
and all the other sets of morphisms between their shifts vanish.
Thus $\tilde C$ is not equivalent to $C$ or $C_1$.
However, we can prove that $\tilde C$ is indeed not admissible:

\begin{Lem}\label{not-saturated}
$\tilde C$ is not saturated.
\end{Lem}

\begin{proof}
Since $\mathcal{O}_E(-1,0)$ is exceptional, it generates a saturated 
subcategory of $\tilde C$.
Therefore, it is sufficient to prove that its left orthogonal
\[
\tilde C' = {}^{\perp}\langle \mathcal{O}_E(-1,0) \rangle \subset \tilde C
\]
is not saturated.
We decompose the other object $\mathcal{O}_E(0,-1)$ 
by a distinguished triangle
\[
c_1 \to \mathcal{O}_E(0,-1) \to \mathcal{O}_E(-1,0)[2] \to c_1[1].
\]
where $c_1 \in \tilde C'$.
Then $\tilde C'$ is generated by $c_1$, and we have
\[
\text{Hom}(c_1,c_1[p]) \cong \begin{cases} k \quad \text{ if } p = 0,3 \\
0 \quad \text{ otherwise.} \end{cases}
\]
Thus $\tilde C'$ is equivalent to a category generated by an object
$\mathcal{O}_M$ in the derived category $D^b(\text{Coh}(M))$ 
for a Calabi-Yau $3$-fold $M$.
Such an object is called a {\em spherical object}.
Since $D^b(\text{Coh}(M))$ is indecomposable and its Serre functor is 
isomorphic to a shift functor $[3]$, 
$\langle \mathcal{O}_M \rangle$ is not saturated,
hence neither is $\tilde C'$.

We can also check this fact directly.
Let $F: (\tilde C')^{\text{op}} \to (k\text{-vect})$ be a cohomological functor
such that $F(c_1) \cong k$ and $F(c_1[p]) = 0$ for $p \ne 0$.
Then $F$ is not representable by any object $c_2 \in \tilde C'$.
Indeed, if $F(c) \cong \text{Hom}(c,c_2)$ for arbitrary $c \in \tilde C'$, then
$\sum_p (-1)^p \dim \text{Hom}(c_1[p],c_2)$ should be even, a contradiction.
\end{proof}


\subsection{Case 2}

$E$ is isomorphic to a projective plane  
$\mathbb{P}^2$, and the normal bundle of $E$ is isomorphic
to $\mathcal{O}_E(-2)$.

Let 
\[
C = \langle \mathcal{O}_E(-1) \rangle \subset D^b(\text{Coh}(Y)).
\]
Since $\mathcal{O}_E(-1)$ is an exceptional object, $C$ is an 
admissible subcategory.
We note that $(\mathcal{O}_E(-2), \mathcal{O}_E(-1))$ is not an 
exceptional collection as in Case 1.
Indeed, we have
\[
\text{Hom}(\mathcal{O}_E(-1), \mathcal{O}_E(-2)[3]) 
\cong \text{Hom}(\mathcal{O}_E(-2), \mathcal{O}_E(-2))^* 
\ne 0.
\]
Thus these two objects may not generate an admissible subcategory.

Let $\mathcal{X}$ be the smooth Deligne-Mumford stack associated to
the variety $X$ which has only a quotient singularity $P$, and 
let $\pi: \mathcal{X} \to X$ be the projection.
We know that $D \cong D^b(\text{Coh}(\mathcal{X}))$ 
by \cite{toric}, and the latter is the correct category.
We shall identify these categories in the following.

We have $f_*\mathcal{O}_P(1) = 0$.
But $\mathcal{O}_P(1)$ is not an exceptional object, because
$\text{Hom}^2(\mathcal{O}_P(1),\mathcal{O}_P(1)) \ne 0$.
We know also that $S_D^2(d) \cong d[6]$ if $d \in D$ and $f_*d = 0$, 
but 
\[
S_D(\mathcal{O}_P(1)) \cong \mathcal{O}_P[3] 
\not\cong \mathcal{O}_P(1)[3].
\] 

\begin{Lem}
$D$ is minimal.
\end{Lem}

\begin{proof}
Suppose that there is a semi-orthgonal decomposition
$D = \langle C',D' \rangle$ such that $f^*\text{Perf}(X) \subset D'$.
Let $c \in C'$ be a non-zero object.
Since $f_*c=0$, $c$ is supported at the point $P$.
Since $\text{Hom}(\mathcal{O}_{\mathcal{X}}, c[n])=0$ for any $n$,
$H^n(c)$ should be of the form $\mathcal{O}_P(1)^{c_n}$.
Any object $d \in D'$ has a finite resolution whose terms are 
of the form $\mathcal{O}_{\mathcal{X}}^{a_n} \oplus 
\mathcal{O}_{\mathcal{X}}(1)^{b_n}$ near $P$.
If $b_n \ne 0$ for some $n$ and for any choice of such a resolution, 
then we have $\text{Hom}(d,c[m]) \ne 0$ for some $m$, a contradiction.
Therefore, $b_n = 0$ for all $n$, 
and $d \in f^*\text{Perf}(X)$, a contradiction to Lemma~\ref{saturated}.
\end{proof}


\subsection{Case 3}

$E$ is isomorphic to a singular quadric surface, 
and the normal bundle of $E$ is isomorphic to $\mathcal{O}_E(-1)$.

This case has only a partial answer.
Let $l$ be a ruling. 
Then $\mathcal{O}_E(1) \cong \mathcal{O}_E(2l)$.
There is an exact sequence
\[
0 \to \mathcal{O}_E(-3l) \to \mathcal{O}_E(-2l)^2 \to \mathcal{O}_E(-l)
\to 0.
\]
The Serre functor of $D^b(\text{Coh}(Y))$ is given by 
$\otimes \mathcal{O}_Y(E)[3]$, and we have
\[
\text{Hom}(\mathcal{O}_E(-l), \mathcal{O}_E(-l)[2])
\cong \text{Hom}(\mathcal{O}_E(-l), \mathcal{O}_E(-3l)[1])^*
\ne 0
\]
hence $\mathcal{O}_E(-l)$ is not an exceptional object.

Let 
\[
C = \langle \mathcal{O}_E(-1) \rangle
\]
and $D = {}^{\perp}C$.
Since $\mathcal{O}_E(-1)$ is an exceptional object, $C$ is admissible.
We have $\mathcal{O}_E(-l) \in D$.
The above sequence implies that 
\[
j^!\mathcal{O}_E(-3l) \cong \mathcal{O}_E(-l)[-1]
\] 
where $j^!$ is the right adjoint to $j_*: D \to D^b(\text{Coh}(Y))$.
Therefore, 
\[
S_D(\mathcal{O}_E(-l)) \cong \mathcal{O}_E(-l)[2].
\]

We expect that if $c \in D$ is right orthogonal to $f^*\text{Perf}(X)$,
then $S_D(c) \cong c[2]$.
We note that the shift number $2$ is different from usual $3 = \dim X$.
We can only prove a partial result:

\begin{Lem}
Let $c \in D$ such that $c \in (f^*\text{Perf}(X))^{\perp}$.
Assume in addition that there exists $d \in D^b(\text{Coh}(E))$
such that $c = i_*d$ for $i: E \to Y$.
Then $S_D(c) \cong c[2]$.
\end{Lem}

\begin{proof}
By a generalization of \cite{Beilinson}, $d$ has a two sided resolution
whose terms are direct sums of the sheaves $\mathcal{O}_E$, 
$\mathcal{O}_E(-l)$ and $\mathcal{O}_E(-2l)$.
Since 
\[
\text{Hom}_E(\mathcal{O}_E, d[n]) 
\cong \text{Hom}_Y(\mathcal{O}_Y, c[n])
= 0
\]
for any $n$, the terms do not contain $\mathcal{O}_E$.
On the other hand, since 
\[
\text{Hom}_Y(c, \mathcal{O}_E(-2l)[n]) = 0
\]
for any $n$, the terms do not contain $\mathcal{O}_E(-2l)$ either.
Since $S_D(\mathcal{O}_E(-l)) \cong \mathcal{O}_E(-l)[2]$, we have our
assertion.
\end{proof}

Let $B$ be the full subcategory of $D^b(\text{Coh}(Y))$ 
consisting of objects
whose supports are contained in $E$.
Then $B$ is generated by $\mathcal{O}_E$-modules 
as a triangulated category.
But we note that the inclusion functor $D^b(\text{Coh}(E)) \to B$ 
is not fully faithful 
because there are more extensions in $B$.

In order to prove the above expectation, one would need more geometric 
argument.
Anyway, if it is true, then the minimality 
of the category $D$ will follow from its indecomposability.

\begin{Lem}
$D$ is indecomposable.
\end{Lem}

\begin{proof}
Assume that $D = C' \oplus D'$.
Since $\text{Hom}(\mathcal{O}_Y, \mathcal{O}_Y) \cong k$,
$\mathcal{O}_Y \in D$ is indecomposable.
We may assume that $\mathcal{O}_Y \in D'$.

Let $x \in E$ be an arbitrary point, and let $d = j^!\mathcal{O}_x$ 
for the right adjoint functor $j^!$ of $j: D \to A = D^b(\text{Coh}(Y))$.
We have a distinguished triangle
\[
d \to \mathcal{O}_x \to 
\text{Hom}^{\bullet}(\mathcal{O}_x, \mathcal{O}_E(-1))^* 
\otimes \mathcal{O}_E(-1) \to d[1].
\]
The third term is isomorphic to $\mathcal{O}_E(-1)[2] \oplus 
\mathcal{O}_E(-1)[3]$.
Then
\[
\text{Hom}(d,d) \cong \text{Hom}(d,\mathcal{O}_x) 
\cong \text{Hom}(\mathcal{O}_x,\mathcal{O}_x) \cong k
\]
because 
\[
\text{Hom}(\mathcal{O}_E(-1)[k],\mathcal{O}_x) = 0
\]
for $k > 0$.
Therefore, $d$ is indecomposable.
Since 
\[
\text{Hom}(\mathcal{O}_Y, d) 
\cong \text{Hom}(\mathcal{O}_Y,\mathcal{O}_x) 
\ne 0
\]
we have $d \in D'$.

Let $c \in C'$ be an arbitrary object such that $c \not\cong 0$.
Then the support of $c$ is contained in $E$.
Thus there exists $x$ and $p$ such that 
\[
\text{Hom}(c, d[p]) \cong 
\text{Hom}(c,\mathcal{O}_x[p]) \ne 0
\]
a contradiction.
\end{proof}

\section{Questions}

Based on the above observation, 
we would like to ask the following questions: 

(1) Let $X$ be a projective variety with only terminal singularities, 
and $f: Y \to X$ a resolution of singularities.
Does there exist a minimal saturated subcategory of 
$D^b(\text{Coh}(Y))$ which contains $f^*\text{Perf}(X)$?
Are two such minimal saturated subcategories equivalent?

(2) More generally, let $X$ be a smooth projective variety.
Does the category $D^b(\text{Coh}(X))$ have finite length 
with respect to semi-orthogonal decompositions? 
Does a Jordan-H\"older type theorem hold for semi-orthogonal 
decompositions of saturated categories?

(3) It would be nice if we have more method for testing 
the admissibility of a subcategory in general situation 
(cf. Lemma~\ref{not-saturated}).

(4) Can a geometric argument be given to establish the expectation from
Subsection 4.3? (added by a referee)

Department of Mathematical Sciences, University of Tokyo, 

Komaba, Meguro, Tokyo, 153-8914, Japan 

kawamata@ms.u-tokyo.ac.jp

\end{document}